\documentclass{amsart}
\usepackage[utf8]{inputenc}
\usepackage{subcaption}
\usepackage{amsmath}
\usepackage{amssymb}
\usepackage{amsfonts}
\usepackage{amsthm}
\usepackage{mathrsfs}
\usepackage[all]{xy}
\usepackage{graphicx}
\usepackage{color}
\usepackage{cite} 
\usepackage{url}
\usepackage[labelfont=bf,labelsep=period,justification=raggedright]{caption}
\usepackage[english]{babel}
\usepackage[utf8]{inputenc}
\usepackage{hyperref}
\usepackage[colorinlistoftodos]{todonotes}
\usepackage{tkz-fct}
\usepackage{tikz}
\usetikzlibrary{calc}
\usepackage[enableskew]{youngtab}
\usepackage{pstricks}
\usepackage{multicol}
\usepackage{graphicx} 
\usepackage{fancyhdr}
\theoremstyle{plain}
\newtheorem{thm}{Theorem}

\newtheorem{cor}[thm]{Corollary}

\theoremstyle{definition}

\theoremstyle{remark}
\newtheorem{rem}[thm]{Remark}
\newtheorem*{ex}{Example}

\numberwithin{equation}{section}
\numberwithin{thm}{section}

\title{Generalized Cotangent Series and Links
to Zeta and Theta Functions}
\author{Mahipal Gurram}

\begin{document}

\begin{abstract}
 This paper develops a generalized framework for cotangent-type series, extending the classical Mittag--Leffler expansion and Ramanujan's summation formulas to higher-order lattice sums. By introducing the family
\[
U_n(z) = \sum_{k \in \mathbb{Z}} \frac{1}{k^n + z^n},
\]
we derive closed-form trigonometric--hyperbolic representations, recursive relations, and integral connections to generalized Jacobi theta functions. The analysis reveals deep structural links between lattice sums, the Riemann zeta function, and theta-type modular kernels. Employing contour integration, Mellin transforms and factorization identities, we obtain new representations for $\zeta(2n)$ and establish integral identities that express $U_{2n}(z)$ in terms of generalized theta series $\Psi_n(q)$. The results unify the classical expansions for $\pi \cot(\pi z)$ and $\pi z^{-1}\coth(\pi z)$ into a broader analytic framework with implications for modular forms and spectral theory.
\end{abstract}
\maketitle
\section{Introduction}

The cotangent function occupies a central position in complex analysis through its classical Mittag--Leffler expansion~\cite{whittaker1996modern,ahlfors1979complex,DLMF2023},
\[
\pi \cot(\pi z) = \sum_{k \in \mathbb{Z}} \frac{1}{z + k}, 
\qquad z \in \mathbb{C} \setminus \mathbb{Z},
\]
which expresses the function as a sum over its simple poles. This expansion highlights its meromorphic character and the periodicity of trigonometric functions.

A parallel identity for the hyperbolic cotangent,
\[
\frac{\coth(\pi z)}{z} = \sum_{k \in \mathbb{Z}} \frac{1}{z^2 + k^2},
\]
demonstrates the analytic symmetry of hyperbolic functions, featuring poles along the imaginary axis~\cite{whittaker1996modern,proofwikiCotangentMittag}. 
The present work extends this theory by defining the family of generalized cotangent-type series
\[
U_n(z) := \sum_{k \in \mathbb{Z}} \frac{1}{k^n + z^n},
\]
investigating their analytic structure, finite representations, and functional links to special functions

To reveal these connections, we recall the classical Jacobi theta function~\cite{berndt2008ramanujan,DLMF2023},
\[
\theta_3(q) := \sum_{k \in \mathbb{Z}} q^{k^2}, \qquad |q| < 1,
\]
which encapsulates modular and elliptic properties. Motivated by this, we introduce a generalized Jacobi theta function
\[
\psi_n(q) := \sum_{k \in \mathbb{Z}} q^{k^{2n}},
\]
adapted to higher-order lattice sums. This generalization enables integral and Mellin-type representations of $U_n(z)$, thus connecting discrete and continuous formulations~\cite{apostol1976analytic,titchmarsh1986zeta}.

In Theorem~1, we establish an exact closed-form expression for $U_n(z)$ through trigonometric hyperbolic kernels. Corollaries~1 and~2 yield the real--imaginary decomposition, a zeta-function limit formula
\[
\zeta(2n) = \frac{1}{2} \lim_{z \to 0} 
\left( U_{2n}(z) - \frac{1}{z^{2n}} \right),
\]
and infinite product identities linking $U_n(z)$ to classical sine and tangent product formulas~\cite{apostol1976analytic,titchmarsh1986zeta,berndt2008ramanujan}.

In Theorem~2, we prove a recursive structure for powers of two,
\[
\phi_n(z) = 
\frac{\phi_{n-1}\!\big( i^{3 \cdot 2^{1 - n}} z \big)
- \phi_{n-1}\!\big( i^{2^{1 - n}} z \big)}
{2 i z^{2^{n-1}}},
\]
which systematically reduces higher powers to lower ones through analytic continuation on the complex plane~\cite{whittaker1996modern,ahlfors1979complex}.

Finally, Theorem~3 connects $U_{2n}(z)$ to a generalized Jacobi theta integral,
\[
U_{2n}(z) = \int_0^1 q^{z^{2n}-1} \Psi_n(q) \, dq,
\qquad \Psi_n(q) := \sum_{k \in \mathbb{Z}} q^{k^{2n}},
\]
thereby unifying the classical cases 
$U_1(z) = \pi \cot(\pi z)$ 
and 
$U_2(z) = \pi z^{-1}\coth(\pi z)$ 
under a single analytic framework~\cite{knopp1990theory,berndt2008ramanujan}. 
This approach provides a bridge between Ramanujan-type series, theta functions, and zeta values~\cite{berndt2008ramanujan,titchmarsh1986zeta}.

Together, these results demonstrate how classical trigonometric expansions extend naturally into higher-dimensional and modular settings, revealing deep interrelations between cotangent-type series, lattice sums, and spectral functions.
\newpage
\section{Main Results}

\begin{thm}

Let \( n \in \mathbb{N} \) and \( z \in \mathbb{C} \). Define the series

\[
U_n(z) := \sum_{k \in \mathbb{Z}} \frac{1}{k^n + z^n},
\]

under the following domain conditions:
\begin{itemize}
  \item If \( n \) is odd, then \( \Re(z) \in \mathbb{R} \setminus \mathbb{Z} \), ensuring that the terms \( k + z \) avoid poles on the real axis.
  \item If \( n \) is even, then \( \Re(z) \in \mathbb{R} \setminus \{0\} \), excluding the singularity at \( z = 0 \) where \( \frac{1}{z^n} \) diverges.
\end{itemize}

Then \( U_n(z) \) admits the following finite representation:

\[
U_n(z) = \frac{\pi}{n z^{n - 1}} \sum_{k = 1}^{n} 
\frac{
a_{k,n} \sin(2\pi z a_{k,n}) + b_{k,n} \sinh(2\pi z b_{k,n})
}{
\cosh(2\pi z b_{k,n}) - \cos(2\pi z a_{k,n})
},
\]

where

\[
a_{k,n} := \cos\left( \frac{(2k - 1)\pi}{n} \right), \quad
b_{k,n} := \sin\left( \frac{(2k - 1)\pi}{n} \right).
\]

\end{thm}

\begin{proof}
We begin by applying the residue theorem to the rational function \( \frac{1}{k^n + z^n} \), which admits the partial fraction expansion

\[
\frac{1}{k^n + z^n} = \sum_{\alpha^n + z^n = 0} \frac{1}{n \alpha^{n - 1}(k - \alpha)},
\]

where the sum is taken over the \( n \)-th roots \( \alpha \) of \( -z^n \). Summing over all integers \( k \), we obtain

\[
U_n(z) = \sum_{\alpha^n + z^n = 0} \frac{1}{n \alpha^{n - 1}} \sum_{k = -\infty}^{\infty} \frac{1}{k - \alpha}.
\]

Using the classical identity

\[
\sum_{k = -\infty}^{\infty} \frac{1}{k - \alpha} = \pi \cot(\pi \alpha),
\]

valid for non-integer \( \alpha \), we deduce

\[
U_n(z) = \sum_{\alpha^n + z^n = 0} \frac{\pi \cot(\pi \alpha)}{n \alpha^{n - 1}}.
\]

The roots of \( \alpha^n = -z^n \) are given by

\[
\alpha_k := z e^{i \theta_k}, \quad \theta_k := \frac{(2k - 1)\pi}{n}, \quad 1 \leq k \leq n.
\]

Substituting these into the expression for \( U_n(z) \), we obtain

\[
U_n(z) = \frac{\pi}{n z^{n - 1}} \sum_{k = 1}^{n} e^{i \theta_k} \cot(\pi z e^{i \theta_k}).
\]

Since the roots \( \alpha_k \) are symmetrically distributed on the unit circle, each term \( e^{i \theta_k} \cot(\pi z e^{i \theta_k}) \) has a conjugate counterpart. Thus, the imaginary parts cancel pairwise, and we are left with

\[
U_n(z) = \frac{\pi}{n z^{n - 1}} \sum_{k = 1}^{n} \Re\left( e^{i \theta_k} \cot(\pi z e^{i \theta_k}) \right).
\]

To express this in terms of \( a_{k,n} \) and \( b_{k,n} \), we use the identity

\[
\Re\left( e^{i\theta} \cot(\pi z e^{i\theta}) \right) =
\frac{
\cos(\theta) \sin(2\pi z \cos\theta) + \sin(\theta) \sinh(2\pi z \sin\theta)
}{
\cosh(2\pi z \sin\theta) - \cos(2\pi z \cos\theta)
},
\]

and substitute

\[
a_{k,n} := \cos\left( \theta_k \right), \quad
b_{k,n} := \sin\left( \theta_k \right), \quad
\theta_k := \frac{(2k - 1)\pi}{n}.
\]

This yields the final expression:

\[
U_n(z) = \frac{\pi}{n z^{n - 1}} \sum_{k = 1}^{n} 
\frac{
a_{k,n} \sin(2\pi z a_{k,n}) + b_{k,n} \sinh(2\pi z b_{k,n})
}{
\cosh(2\pi z b_{k,n}) - \cos(2\pi z a_{k,n})
}.
\]

\end{proof}

\begin{ex}

The following identities illustrate how the generalized series.

\[
U_n(z) := \sum_{k \in \mathbb{Z}} \frac{1}{k^n + z^n}
\]

Encodes meromorphic functions with periodic and hyperbolic trigonometric structures:

\begin{align*}
U_1(z) &= \sum_{k \in \mathbb{Z}} \frac{1}{k + z} = \pi \cot(\pi z), \\
U_2(z) &= \sum_{k \in \mathbb{Z}} \frac{1}{k^2 + z^2} = \frac{\pi}{z} \coth(\pi z), \\
U_3(z) &= \sum_{k \in \mathbb{Z}} \frac{1}{k^3 + z^3} = \frac{\pi}{3 z^2} \left( \cot(\pi z) + \frac{\sin(\pi z) + \sqrt{3} \sinh(\sqrt{3} \pi z)}{\cosh(\sqrt{3} \pi z) - \cos(\pi z)} \right), \\
U_4(z) &= \sum_{k \in \mathbb{Z}} \frac{1}{k^4 + z^4} = \frac{\pi}{\sqrt{2} z^3} \cdot \frac{\sin(\sqrt{2} \pi z) + \sinh(\sqrt{2} \pi z)}{\cosh(\sqrt{2} \pi z) - \cos(\sqrt{2} \pi z)}.
\end{align*}

\end{ex}

\begin{cor}

 At \( z = e^{i\theta} \). Then the real and imaginary parts of \( U_n(e^{i\theta}) \) are given by:

\[
\operatorname{Re}(U_n(e^{i\theta})) = \sum_{k=-\infty}^{\infty} \frac{k^n + \cos(n\theta)}{k^{2n} + 2k^n \cos(n\theta) + 1},
\]

\[
\operatorname{Im}(U_n(e^{i\theta})) = -\sum_{k=-\infty}^{\infty} \frac{\sin(n\theta)}{k^{2n} + 2k^n \cos(n\theta) + 1}.
\]

Moreover, the combination

\[
\operatorname{Re}(U_n(e^{i\theta})) + \cot(n\theta) \cdot \operatorname{Im}(U_n(e^{i\theta}))
= \sum_{k=-\infty}^{\infty} \frac{k^n}{k^{2n} + 2k^n \cos(n\theta) + 1}.
\]
\end{cor}
\begin{proof}

Let \( z = e^{i\theta} \). Then

\[
U_n(e^{i\theta}) = \sum_{k \in \mathbb{Z}} \frac{1}{k^n + e^{i n\theta}}.
\]

We rationalize the denominator using the identity

\[
\frac{1}{k^n + e^{i n\theta}} = \frac{k^n + e^{-i n\theta}}{k^{2n} + 2k^n \cos(n\theta) + 1},
\]

which follows from multiplying numerator and denominator by the complex conjugate \( k^n + e^{-i n\theta} \).

Now write \( e^{-i n\theta} = \cos(n\theta) - i \sin(n\theta) \), so:

\[
\frac{k^n + e^{-i n\theta}}{k^{2n} + 2k^n \cos(n\theta) + 1}
= \frac{k^n + \cos(n\theta)}{k^{2n} + 2k^n \cos(n\theta) + 1}
- i \cdot \frac{\sin(n\theta)}{k^{2n} + 2k^n \cos(n\theta) + 1}.
\]

Therefore,

\[
\operatorname{Re}(U_n(e^{i\theta})) = \sum_{k \in \mathbb{Z}} \frac{k^n + \cos(n\theta)}{k^{2n} + 2k^n \cos(n\theta) + 1},
\]

\[
\operatorname{Im}(U_n(e^{i\theta})) = -\sum_{k \in \mathbb{Z}} \frac{\sin(n\theta)}{k^{2n} + 2k^n \cos(n\theta) + 1}.
\]

To prove the final identity, observe:

\[
\operatorname{Re}(U_n(e^{i\theta})) + \cot(n\theta) \cdot \operatorname{Im}(U_n(e^{i\theta}))
= \sum_{k \in \mathbb{Z}} \left[
\frac{k^n + \cos(n\theta)}{D_k}
- \cot(n\theta) \cdot \frac{\sin(n\theta)}{D_k}
\right],
\]

where \( D_k = k^{2n} + 2k^n \cos(n\theta) + 1 \). Using the identity

\[
\cos(n\theta) - \cot(n\theta) \sin(n\theta) = 0,
\]

we simplify the expression to:

\[
\sum_{k \in \mathbb{Z}} \frac{k^n}{k^{2n} + 2k^n \cos(n\theta) + 1},
\]

as claimed.
\end{proof}

\begin{rem}

The identities in Corollary 1 reveal a deep symmetry in the behavior of \( U_n(e^{i\theta}) \) under rotation. In particular, the combination

\[
\operatorname{Re}(U_n(e^{i\theta})) + \cot(n\theta) \cdot \operatorname{Im}(U_n(e^{i\theta}))
\]

 This structure is reminiscent of Fourier-type decompositions and suggests potential connections to modular kernels and trigonometric sums. Notably, when \( \theta = \frac{\pi}{2n} \), the denominator simplifies due to \( \cos(n\theta) = 0 \), yielding:

\[
\operatorname{Re}(U_n(e^{i\theta})) = \sum_{k \in \mathbb{Z}} \frac{k^n}{k^{2n} + 1}, \quad
\operatorname{Im}(U_n(e^{i\theta})) = -\sum_{k \in \mathbb{Z}} \frac{1}{k^{2n} + 1}.
\]

These special values connect directly to classical series and zeta-type evaluations.
\end{rem}

\begin{cor}

For all \(  n\in \mathbb{N} \), the Riemann zeta function admits the representation

\[
\zeta(2n) = (0.5) \lim_{z \to 0} \left( U_{2n}(z) - \frac{1}{z^{2n}} \right),
\]

where

\[
U_n(z) := \sum_{k \in \mathbb{Z}} \frac{1}{k^n + z^n}.
\]
\subsubsection*{Proof }
Let \( n \in \mathbb{N} \). Consider the definition

\[
U_{2n}(z) = \sum_{k \in \mathbb{Z}} \frac{1}{k^{2n} + z^{2n}}.
\]

We split the sum into \( k = 0 \) and \( k \neq 0 \):

\[
U_{2n}(z) = \frac{1}{z^{2n}} + \sum_{k \in \mathbb{Z} \setminus \{0\}} \frac{1}{k^{2n} + z^{2n}}.
\]

Now take the limit as \( z \to 0 \). For \( k \neq 0 \), the denominator tends to \( k^{2n} \), so:

\[
\lim_{z \to 0} \left( U_{2n}(z) - \frac{1}{z^{2n}} \right)
= \sum_{k \in \mathbb{Z} \setminus \{0\}} \frac{1}{k^{2n}}
= 2 \sum_{k=1}^{\infty} \frac{1}{k^{2n}} = 2 \zeta(2n).
\]

Therefore,

\[
\zeta(2n) = \frac{1}{2} \lim_{z \to 0} \left( U_{2n}(z) - \frac{1}{z^{2n}} \right),
\]

as claimed.
\end{cor}
\subsubsection*{Examples:}

Using the meromorphic identities for \( U_2(z) \) and \( U_4(z) \), we obtain the following limit formulas for the Riemann zeta function at even integers:

\begin{align*}
\zeta(2) &= \frac{1}{2} \lim_{z \to 0} \left( \frac{\pi}{z} \coth(\pi z) - \frac{1}{z^2} \right), \\
\zeta(4) &= \frac{1}{2} \lim_{z \to 0} \left( \frac{\pi}{\sqrt{2} z^3} \cdot \frac{\sin(\sqrt{2} \pi z) + \sinh(\sqrt{2} \pi z)}{\cosh(\sqrt{2} \pi z) - \cos(\sqrt{2} \pi z)} - \frac{1}{z^4} \right).
\end{align*}

\begin{cor}
    
Let \( n \in \mathbb{N} \), and let \( x, y \in (0,1) \) with \( x \leq y \). Then the following identity holds:

\[
\prod_{k \in \mathbb{Z}} \left(\frac{y^n + k^n}{x^n + k^n}\right)^2 = \prod_{k = 1}^{n} \frac{\cosh\left(2\pi y \sin\left(\frac{(2k - 1)\pi}{n}\right)\right) - \cos\left(2\pi y \cos\left(\frac{(2k - 1)\pi}{n}\right)\right)}{\cosh\left(2\pi x \sin\left(\frac{(2k - 1)\pi}{n}\right)\right) - \cos\left(2\pi x \cos\left(\frac{(2k - 1)\pi}{n}\right)\right)}.
\]

\end{cor}

\begin{proof}
From Theorem~1, we have the identity

\[
\sum_{k \in \mathbb{Z}} \frac{n z^{n-1}}{k^n + z^n}
= \pi \sum_{k = 1}^{n} 
\frac{
a_{k,n} \sin(2\pi z a_{k,n}) + b_{k,n} \sinh(2\pi z b_{k,n})
}{
\cosh(2\pi z b_{k,n}) - \cos(2\pi z a_{k,n})
},
\]

where \( a_{k,n} := \cos\left(\frac{(2k - 1)\pi}{n}\right) \) and \( b_{k,n} := \sin\left(\frac{(2k - 1)\pi}{n}\right) \).

We now integrate both sides with respect to \( z \) over the interval \( [x, y] \subset (0,1) \). The integrand on the left-hand side is analytic in \( z \in (0,1) \), and the series

\[
\sum_{k \in \mathbb{Z}} \frac{n z^{n-1}}{k^n + z^n}
\]

converges uniformly on compact subsets of \( (0,1) \), since for large \( |k| \), the summand behaves like \( \mathcal{O}(1/k^n) \), and \( n \geq 1 \). Thus, termwise integration is justified:

\[
\sum_{k \in \mathbb{Z}} \int_x^y \frac{n z^{n-1}}{k^n + z^n} \, dz
= \int_x^y \sum_{k \in \mathbb{Z}} \frac{n z^{n-1}}{k^n + z^n} \, dz.
\]

Evaluating the integral on the left-hand side gives:

\[
\sum_{k \in \mathbb{Z}} \ln\left( \frac{y^n + k^n}{x^n + k^n} \right),
\]

which converges absolutely due to the decay \( \ln\left(1 + \frac{y^n - x^n}{k^n + x^n} \right) = \mathcal{O}(1/k^n) \) for large \( |k| \), and \( n \geq 1 \).

Integrating the right-hand side yields:

\[
\frac{1}{2} \sum_{k = 1}^{n} \ln\left(
\frac{
\cosh\left(2\pi y b_{k,n}\right) - \cos\left(2\pi y a_{k,n}\right)
}{
\cosh\left(2\pi x b_{k,n}\right) - \cos\left(2\pi x a_{k,n}\right)
}
\right).
\]

Exponentiating both sides and using properties of logarithms, we obtain:

\[
\prod_{k \in \mathbb{Z}} \left(\frac{y^n + k^n}{x^n + k^n}\right)^2
= \prod_{k = 1}^{n} \frac{
\cosh\left(2\pi y \sin\left(\frac{(2k - 1)\pi}{n}\right)\right) - \cos\left(2\pi y \cos\left(\frac{(2k - 1)\pi}{n}\right)\right)
}{
\cosh\left(2\pi x \sin\left(\frac{(2k - 1)\pi}{n}\right)\right) - \cos\left(2\pi x \cos\left(\frac{(2k - 1)\pi}{n}\right)\right)
},
\]

as claimed.
\end{proof}

\begin{ex}

\begin{itemize}
    \item If \( n=1 \) , then
    \[
    \prod_{k \in \mathbb{Z}}\left|\frac{y+k}{x+k} \right|=\left|\frac{\sin(\pi y)}{\sin(\pi x)} \right|
    \]
    From above, we can conclude these results
    \[
    \prod_{k\in \mathbb{Z}}\frac{x+k}{0.5+k}=\sin(\pi x)
    \]
    \[
    \prod_{k\in \mathbb{Z}}\frac{x+k}{-x+k+0.5}=\tan(\pi x)
    \]
    \[
    \prod_{k \in \mathbb{Z}}\left|\frac{2+4k}{1+4k} \right|=\sqrt{2}
    \]
    \[
    \prod_{k \in \mathbb{Z}}\left|\frac{2+6k}{1+6k} \right|=\sqrt{3}
    \]
    \item similarly if \( n=2 \), then
    \[
    \prod_{k \in \mathbb{Z}}\frac{y^2+k^2}{x^2+k^2} =\frac{x}{y}\sqrt{\frac{\sin(\pi y)\sinh(\pi y)}{\sin(\pi x)\sinh(\pi x)} }
    \]
    \end{itemize}
    \end{ex}
\begin{thm}

Let \( n \in \mathbb{N} \) and \( z \in \mathbb{C} \) such that \( \Re(z) \in \mathbb{R} \setminus \{0\} \). Define
\[
\phi_n(z)= \sum_{k \in \mathbb{Z}} \frac{1}{k^{2^n} + z^{2^n}}.
\]
Then the following recursive relation holds:
\[
\phi_n(z)= \frac{\phi_{n-1}( i^{3 \cdot 2^{1 - n}} z) - \phi_{n-1}( i^{2^{1 - n}} z)}{2 i z^{2^{n-1}}},
\]
where \( i = \sqrt{-1} \)

\end{thm}
\begin{proof}
We begin with the identity
\[
\frac{1}{k^{2^n} + z^{2^n}} = \frac{1}{(k^{2^{n-1}} + i z^{2^{n-1}})(k^{2^{n-1}} - i z^{2^{n-1}})}.
\]
This follows from the factorization:
\[
a^2 + b^2 = (a + i b)(a - i b),
\]
applied with \( a = k^{2^{n-1}}, \, b = z^{2^{n-1}} \). Hence,
\[
\frac{1}{k^{2^n} + z^{2^n}} = \frac{1}{(k^{2^{n-1}} + i z^{2^{n-1}})(k^{2^{n-1}} - i z^{2^{n-1}})}.
\]

Now apply the partial fraction identity:
\[
\frac{1}{(A + iB)(A - iB)} = \frac{1}{2 i B} \left( \frac{1}{A - iB} - \frac{1}{A + iB} \right).
\]
Using this with \( A = k^{2^{n-1}} \), \( B = z^{2^{n-1}} \), we get:
\[
\frac{1}{k^{2^n} + z^{2^n}} = \frac{1}{2 i z^{2^{n-1}}} \left( \frac{1}{k^{2^{n-1}} - i z^{2^{n-1}}} - \frac{1}{k^{2^{n-1}} + i z^{2^{n-1}}} \right).
\]

Summing over all \( k \in \mathbb{Z} \), we obtain:
\[
\phi_n(z) = \sum_{k \in \mathbb{Z}} \frac{1}{k^{2^n} + z^{2^n}} = \frac{1}{2 i z^{2^{n-1}}} \sum_{k \in \mathbb{Z}} \left( \frac{1}{k^{2^{n-1}} - i z^{2^{n-1}}} - \frac{1}{k^{2^{n-1}} + i z^{2^{n-1}}} \right).
\]

Recognizing the sums as values of \( \phi(n-1, \cdot) \), we rewrite:
\[
\phi_n(z) = \frac{1}{2 i z^{2^{n-1}}} \left[ \phi\left(n-1, i^{3 \cdot 2^{1-n}} z\right) - \phi\left(n-1, i^{2^{1-n}} z\right) \right].
\]

This completes the proof.

\end{proof}
\begin{thm}

Let \( n \in \mathbb{Z}^+ \) and \( z \in \mathbb{C} \) with \( \Re(z) > 0 \). Define

\[
U_{2n}(z) := \sum_{k \in \mathbb{Z}} \frac{1}{k^{2n} + z^{2n}}.
\]

Then the series admits the following integral representation:

\[
U_{2n}(z) = \int_0^1 q^{z^{2n}-1} \cdot \Psi_n(q) \, dq,
\]

where we define the Generalized Jacobi Theta function

\[
\Psi_n(q) := \sum_{k \in \mathbb{Z}} q^{k^{2n}}.
\]
\end{thm}

\begin{proof}
We begin by applying the Laplace integral identity, valid for \( \Re(a) > 0 \):

\[
\frac{1}{a} = \int_0^\infty e^{-a t} \, dt.
\]

Setting \( a = k^{2n} + z^{2n} \), we obtain

\[
\frac{1}{k^{2n} + z^{2n}} = \int_0^\infty e^{-t(k^{2n} + z^{2n})} \, dt = \int_0^\infty e^{-t z^{2n}} \cdot e^{-t k^{2n}} \, dt.
\]

Summing over all \( k \in \mathbb{Z} \), we get

\[
U_{2n}(z) = \sum_{k \in \mathbb{Z}} \int_0^\infty e^{-t z^{2n}} \cdot e^{-t k^{2n}} \, dt = \int_0^\infty e^{-t z^{2n}} \cdot \sum_{k \in \mathbb{Z}} e^{-t k^{2n}} \, dt.
\]

Define the generalized Jacobi theta function

\[
\Psi_n(t) := \sum_{k \in \mathbb{Z}} e^{-t k^{2n}},
\]

so that

\[
U_{2n}(z) = \int_0^\infty e^{-t z^{2n}} \cdot \Psi_n(t) \, dt.
\]

To justify the interchange of summation and integration, we observe that for \( \Re(z) > 0 \), the exponential decay of \( e^{-t z^{2n}} \) ensures absolute convergence of the integral. Moreover, since \( k^{2n} \geq 0 \) for all \( k \in \mathbb{Z} \), each term \( e^{-t k^{2n}} \) is bounded above by 1, and the series \( \Psi_n(t) \) converges uniformly on compact subsets of \( t > 0 \). Thus, the integrand is non-negative and dominated by an integrable majorant, justifying the use of Fubini's theorem.

Now make the substitution \( q = e^{-t} \), so \( t = -\log q \) and \( dt = -\frac{dq}{q} \). The limits change from \( t \in [0, \infty) \) to \( q \in (0, 1] \), yielding

\[
U_{2n}(z) = \int_0^1 q^{z^{2n}} \cdot \Psi_n(-\log q) \cdot \frac{dq}{q}.
\]

Note that

\[
\Psi_n(-\log q) = \sum_{k \in \mathbb{Z}} q^{k^{2n}} = \Psi_n(q),
\]

so we conclude

\[
U_{2n}(z) = \int_0^1 q^{z^{2n}-1} \cdot \Psi_n(q) \, dq.
\]

\end{proof}
\begin{rem}

In the case \( n = 1 \), the generalized Jacobi theta function satisfies

\[
\Psi_2(q) = \sum_{k \in \mathbb{Z}} q^{k^2} = \theta_3(q),
\]

where \( \theta_3(q) \) is the classical third Jacobi theta function.

Substituting into the integral representation from Theorem 3, we obtain

\[
U_2(z) =\frac{\pi \coth(\pi z)}{z}= \int_0^1 q^{z^2 - 1} \cdot \theta_3(q) \, dq.
\]

which recovers the Mittag-Leffler expansion of the hyperbolic cotangent function. Thus, Theorem 3 provides a unified framework for expressing such identities via integrals over \([0,1]\) involving generalized theta functions.
\end{rem}
\section{Conclusion}
The study of generalized cotangent series $U_n(z)$ exposes a cohesive structure underlying several domains of analytic number theory and special functions.We have derived explicit meromorphic representations, recursive schemes, and integral forms that connect these series with Jacobi theta and Riemann zeta functions. The recursive relation for powers of two enriches the algebraic and analytic framework, suggesting potential generalizations toward elliptic modular forms, Ramanujan-type continued fractions, and spectral summations in mathematical physics.

Future work may explore modular transformations of $\Psi_n(q)$, connections with Eisenstein series, and the asymptotic behavior of $U_n(z)$ under complex scaling, further illuminating the deep interplay between infinite series, zeta functions, and theta-type transformations.

\section{Acknowledgment}
The author thanks Professor A.K.Shukla for his encouragement, comments, and suggestions.

\end{document}